\documentclass{amsproc}

\newtheorem{theorem}{Theorem}
\newtheorem{lemma}[theorem]{Lemma}
\newtheorem{proposition}[theorem]{Proposition}
\newtheorem{corollary}[theorem]{Corollary}
\theoremstyle{definition}

\theoremstyle{remark}

\numberwithin{equation}{section}

\newcommand{\Z}{\mathbb{Z}}

\newcommand{\uloopr}[1]{\ar@'{@+{[0,0]+(-4,5)}@+{[0,0]+(0,10)}@+{[0,0] +(4,5)}}^{#1}}
\newcommand{\uloopd}[1]{\ar@'{@+{[0,0]+(5,4)}@+{[0,0]+(10,0)}@+{[0,0]+ (5,-4)}}^{#1}}
\newcommand{\dloopr}[1]{\ar@'{@+{[0,0]+(-4,-5)}@+{[0,0]+(0,-10)}@+{[0, 0]+(4,-5)}}_{#1}}
\newcommand{\dloopd}[1]{\ar@'{@+{[0,0]+(-5,4)}@+{[0,0]+(-10,0)}@+{[0,0 ]+(-5,-4)}}_{#1}}
\newcommand{\luloop}[1]{\ar@'{@+{[0,0]+(-8,2)}@+{[0,0]+(-10,10)}@+{[0, 0]+(2,2)}}^{#1}}

\usepackage{amssymb}
\usepackage{amsfonts}
\usepackage{latexsym}
\usepackage{amscd}
\usepackage[mathscr]{euscript}
\usepackage{xy} \xyoption{all}

\begin{document}

\title{Leavitt path algebras of Cayley graphs \\ arising from cyclic groups }
\author{Gene Abrams}
\address{Department of Mathematics, University of Colorado,
Colorado Springs, CO 80918 U.S.A.}
%\curraddr{Department of Mathematics and Statistics,
%Case Western Reserve University, Cleveland, Ohio 43403}
\email{abrams@math.uccs.edu}
\thanks{The first  author is partially supported by a Simons Foundation Collaboration Grants for
Mathematicians Award \#208941.  Part of this work formed the basis of the second author's Master of Science presentation at the University of Colorado Colorado Springs, April 2013.  The authors are grateful to Attila Egri-Nagy, who brought to the attention of the first author the potential connection between Leavitt path algebras and Cayley graphs  during the conference ``Graph C*-algebras, Leavitt path algebras and symbolic dynamics", held at University of Western Sydney, February 2013.}

\author{Benjamin Schoonmaker}
\address{Department of Mathematics, Brigham Young University,
Provo, UT  84602 U.S.A.}
\email{ben.schoonmaker@math.byu.edu}

%    General info
\date{October 15, 2013}



\begin{abstract}
For any positive integer $n$ we describe the Leavitt path algebra of the Cayley graph $C_n$ corresponding to the cyclic group $\Z/n\Z$.  Using a Kirchberg-Phillips-type realization result, we show that there are exactly four isomorphism classes of such Leavitt path algebras, arising as the algebras corresponding to the graphs $C_i$ ($3\leq i \leq 6$). 
\end{abstract}

\maketitle

For each finite group $H$ the {\it Cayley graph} $C_H$ of $H$ is a directed graph which encodes information about the relationships between elements of $H$ and a set of generators of $H$.   In the particular case where $H_n = \Z / n\Z$ (and $n\geq 3$), the Cayley graph $C_{H_n}$ (which we denote simply by $C_n$) consists of $n$ vertices $\{v_1, v_2, \hdots, v_n\}$  and $2n$ edges $\{e_1, e_2, \hdots, e_n, f_1, f_2, \hdots, f_n\}$ for which 
$s(e_i) = v_i,  \  r(e_i) = v_{i+1},  \  s(f_i) = v_i,   \ r(f_i) = v_{i-1},$
where indices are interpreted ${\rm mod } \ n$, and where $s(e )$ (resp., $r(e )$) denotes the source (resp., range) vertex of the edge $e$.    
% (So, for instance, $r(e_n) = e_1$, and $r(f_1) = v_n$.) 
 (More precisely, the graph $C_n$ described here is the Cayley graph for the group $\Z / n\Z $ {\it with respect to the subset}  $\{1, n-1\}$.)  So, for instance,  $C_3$ is the graph
 $$  C_3 = \ \ \  {
\def\labelstyle{\displaystyle}
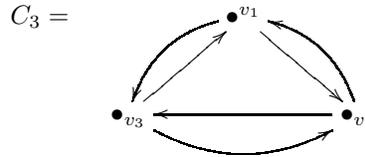
\xymatrix{ {} & \bullet^{v_1}  \ar[rd] \ar@/^{-10pt}/ [ld] &  {} \\
\bullet_{v_3}  \ar[ru] \ar@/^{-15pt}/ [rr]&  & \bullet_{v_2}
 \ar[ll]
\ar@/^{-10pt}/ [lu] \\
}}
$$
\medskip

 \noindent
For the cases $n=1$ and $n=2$, we define the graphs $C_1$ and $C_2$ in a manner consistent with the above description (i.e., as a graph with $n$ vertices and $2n$ edges with appropriate source and range relations), as follows:
$$C_1 = {
\def\labelstyle{\displaystyle}
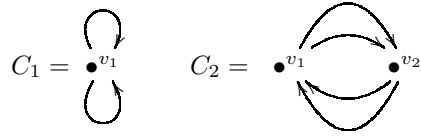
\xymatrix{ \bullet^{v_1}\uloopr{}\dloopr{} 
%& \bullet^{v_2} \ar[l]
}}\ \ \ \ \ \ \  C_2 =   \ \xymatrix{
{\bullet}^{v_1} 
%\ar@(ul,dl)
  \ar@/^1pc/ [r]  \ar@/^2pc/ [r]  & {\bullet}^{v_2} \ar@/^1pc/ [l]  \ar@/^2pc/ [l] }$$

  Less formally, $C_n$ is the graph with $n$ vertices and $2n$ edges, where each vertex emits two edges, one to both of its neighboring vertices.   We denote by $A_{C_n}$ the incidence matrix of $C_n$, which for $n\geq 3$ is easily seen to be the matrix
$$A_{C_n} \ = \ \begin{pmatrix} 0&1&0& \cdots &0&1 \\
        1&0&1 & \cdots & 0&0 \\
         0&1  & 0 & &  0& 0  \\
         & \vdots & &\ddots & \vdots & \\
             0&0& & & 0&1 \\
             1& 0 & & \cdots &1&0 \end{pmatrix}.$$
\noindent
For the cases $n=1$ and $2$ we have $A_{C_1} = (2)$, and $A_{C_2} =  \begin{pmatrix} 0&2 \\
             2& 0  \end{pmatrix}$.  

For any field $K$ and directed graph $E$ the Leavitt path algebra $L_K(E)$ has been the focus of sustained investigation since 2004.  We give here a basic description of $L_K(E)$; for additional information, see e.g.   \cite{AAP1} or \cite{TheBook}.  

\medskip

{\bf Definition of Leavitt path algebra.}   Let $K$ be a field, and let $E = (E^0, E^1, r,s)$ be a directed  graph with vertex set $E^0$ and edge set $E^1$.   The {\em Leavitt path $K$-algebra} $L_K(E)$ {\em of $E$ with coefficients in $K$} is  the $K$-algebra generated by a set $\{v\mid v\in E^0\}$, together with a set of variables $\{e,e^*\mid e\in E^1\}$, which satisfy the following relations:

(V)   \ \ \  \ $vw = \delta_{v,w}v$ for all $v,w\in E^0$, \  
%(i.e., $\{v\mid v\in E^0\}$ is a set of orthogonal idempotents),

  (E1) \ \ \ $s(e)e=er(e)=e$ for all $e\in E^1$,

(E2) \ \ \ $r(e)e^*=e^*s(e)=e^*$ for all $e\in E^1$,

 (CK1) $e^*e'=\delta _{e,e'}r(e)$ for all $e,e'\in E^1$,

(CK2)  $v=\sum _{\{ e\in E^1\mid s(e)=v \}}ee^*$ for every   $v\in E^0$ for which $0 < |s^{-1}(v)| < \infty$.

An alternate description of $L_K(E)$ may be given as follows.  For any graph $E$ let $\widehat{E}$ denote the ``double graph" of $E$, gotten by adding to $E$ an edge $e^*$ for each edge $e\in E^1$.   Then $L_K(E)$ is the usual path algebra $KE^*$, modulo the ideal generated by the relations (CK1) and (CK2).     \hfill $\Box$

\medskip

It is easy to show that $L_K(E)$ is unital if and only if $|E^0|$ is finite.  This is of course the case when $E = C_n$.

Let $E$ be a directed graph  with vertices $v_1, v_2, \hdots, v_n$ and adjacency matrix $A_E = (a_{i,j})$.  We let $F_n$ denote the free abelian monoid on the generators $v_1, v_2, \hdots, v_n$ (so $F_n \cong \oplus_{i=1}^n \Z^+$ as monoids).  We denote the identity element of this monoid by $z$.    We let $R_n(E)$ denote the submonoid of $F_n$ generated by the relations 
$$v_i  = \sum_{j=1}^n a_{i,j}v_j$$
for each non-sink $v_i$.     
%We denote the zero element of $M_E$ by $z$. 
%Equivalently, $M_E$ may be interpreted  as follows.     For each non-sink vertex $v_i$ ($1 \leq i \leq t$)  define  $$\vec{b_i} =  (0, 0, ... , 1, 0, ... 0) \in (\Z^+)^{n}.$$   Consider the equivalence relation  $\sim_E$ in   $(\Z^+)^{n}$,  generated by setting  $$\vec{b_i} \sim_E (a_{i,1}, a_{i,2}, ..., a_{i,n})$$ for each non-sink $v_i$.  Define $$M_E = (\Z^+)^{n} / \sim_E.$$ The operation $+$ in $M_E$ is:  \ \  $[\vec{a}] + [\vec{a'}] = [\vec{a} + \vec{a'}]$   for $\vec{a}, \vec{a'} \in (\Z^+)^{n}$.        We denote the zero element $[(0,0,\hdots,0)]$ of $M_E$  by $z$. 
The {\it graph monoid} $M_E$ of   $E$ is defined as the quotient monoid $$M_E = F_n / R_n(E).$$  The elements of $M_E$ are typically denoted using brackets.   

\medskip

As a representative example, we explicitly describe the graph monoid $M_{C_3}$ associated to the Cayley graph $C_3$.    This is the free abelian monoid on the generators $v_1, v_2, v_3$, modulo the submonoid generated by the relations $v_1 = v_2 + v_3, v_2 = v_1 + v_3, $ and $v_3 = v_1 + v_2$.   Note that, for instance, $v_1 + (v_1 + v_2 + v_3)  = (v_1 + v_2) + (v_1 + v_3) = v_3 + v_2 = v_1$, so that $[v_1] = [v_1] + [v_1 + v_2 + v_3] $ in $M_{C_3}$.  Let $x$ denote the element $[v_1 + v_2 + v_3] $ of $M_{C_3}$.  Then a similar computation yields  that  $[v_2] = [v_2] + x $ and $[v_3] = [v_3] + x$ in $M_{C_3}$.   Moreover, $[v_1] + [v_1] = [v_1] + [v_2 + v_3] = x$, and in a similar fashion we also have $2[v_2] = 2[v_3] = x$ in $M_{C_3}$.   Thus we see that
$$M_{C_3} = \{[z], \ [v_1], \  [v_2], \ [v_3], \ [v_1]+[v_2]+[v_3] \}.$$   (We have not justified why these five elements are distinct in $M_{C_3}$, but this can be done easily;  see e.g. \cite[page 171]{AS}.)

 \medskip
 
More generally, for any $n\geq 1$, the monoid  $M_{C_n}$ is generated by $[v_1], [v_2], \hdots, [v_n]$, subject to the relations $$[v_i ]= [v_{i-1}] + [v_{i+1}]$$ 
(for all $1\leq i \leq n$), where subscripts are interpreted ${\rm mod} \ n$.  (This description also covers the cases $n=1$ and $n=2$.) 

\medskip

  We present now the background information required to achieve our main result (Theorem \ref{MainTheorem}), which yields a description of the Leavitt path algebras corresponding to the Cayley graphs $ \{C_n \ | \ n\geq 1\}$.    The cornerstone of the result is a utilization of the Algebraic Kirchberg Phillips Theorem.   To motivate and explain how this theorem is used, in the following three paragraphs we make a streamlined visit to three elegant, fundamental results in the theory of Leavitt path algebras and purely infinite simple algebras.   
  
  \medskip

For a unital $K$-algebra $A$,   the set of isomorphism classes of finitely generated projective left $A$-modules is denoted by $\mathcal{V}(A)$.  We denote the elements of  $\mathcal{V}(A)$ using brackets; for example, $[A] \in \mathcal{V}(A)$ represents the isomorphism class of the left regular module ${}_AA$.   $\mathcal{V}(A)$ is a monoid, with operation $\oplus$, and zero element $[\{0\}]$.   The monoid $(\mathcal{V}(A), \oplus)$ is {\it conical}; this means that the sum of any two nonzero elements of $\mathcal{V}(A)$ is nonzero, or, rephrased, that $\mathcal{V}(A)^* = \mathcal{V}(A) \setminus \{0\}$ is a semigroup under $\oplus$. 
A striking property of Leavitt path algebras  was established in \cite[Theorem 3.5]{AMP}, to wit:  
\begin{equation}
  \mathcal{V}(L_K(E)) \cong M_E \ \mbox{as monoids}.    \tag{$*$}
 \end{equation}
\vspace{-.2in}
$$\mbox{Moreover, } \   [L_K(E)] \leftrightarrow \sum_{v\in E^0} [v]  \ \mbox{under this isomorphism.}$$

\medskip

A unital $K$-algebra $A$ is called {\it purely infinite simple} in case $A$ is not a division ring, and $A$ has the property that for every nonzero element $x$ of  $A$ there exists $b,c\in A$ for which $bxc=1_A$.     It is shown in \cite[Corollary 2.2]{AGP} that if $A$ is  a unital purely infinite simple $K$-algebra, then the semigroup $(\mathcal{V}(A)^*, \oplus)$ is in fact a group, and, moreover, that $\mathcal{V}(A)^* \cong K_0(A)$, the Grothendieck group of $A$.     (Indeed, for unital Leavitt path algebras, the converse is true as well: if $\mathcal{V}(L_K(E))^*$ is a group, then $L_K(E)$ is purely infinite simple.  This converse is not true for general $K$-algebras.)   Summarizing: when $L_K(E)$ is unital purely infinite simple we have the following isomorphism of groups:  
\begin{equation}
 K_0(L_K(E)) \cong \mathcal{V}(L_K(E))^* \cong M_E^*. \tag{$**$}
 \end{equation}

 \medskip

The finite graphs $E$ for which the Leavitt path algebra $L_K(E)$ is purely infinite simple have been explicitly described in \cite{AAP2}, to wit: 
\begin{equation}
 L_K(E) \  \mbox {is purely infinite simple} \ \Longleftrightarrow  \tag{$***$}
 \end{equation}
\vspace{-.23in}
$$  \ E \  \mbox{is cofinal, sink-free, and satisfies Condition (L)}.$$
 \noindent
Somewhat more fully,  these are the graphs $E$ for which:      every vertex in $E$ connects (via some directed path) to every cycle of $E$;    every cycle in $E$ has an {\it exit} (i.e., in each cycle of $E$ there is a vertex which emits at least two edges); and  $E$ contains at least one cycle.   (The structure of the field $K$ plays no role in determining the purely infinite simplicity of $L_K(E)$.)    

\medskip

We now have the necessary background information in hand which allows us to present the powerful tool which will yield our main result.    The proof of the following theorem utilizes deep results and ideas in the theory of symbolic dynamics. The letters K and P in its name derive from E. Kirchberg and N.C. Phillips, who (independently in 2000) proved an analogous result for graph C$^*$-algebras.  

\medskip

 {\bf The Algebraic KP Theorem.} \cite[Corollary 2.7]{ALPS}  
Suppose $E$ and $F$ are finite graphs for which the Leavitt path algebras $L_K(E)$ and $L_K(F)$ are purely infinite simple.   Suppose that there is an isomorphism $\varphi : K_0(L_K(E)) \rightarrow K_0(L_K(F))$ for which $\varphi([L_K(E)]) = [L_K(F)]$, and suppose also that the two integers ${\rm det}(I_{|E^0|}~-~A_E^t)$ and ${\rm det}~(I_{|F^0|}~-~A_F^t)$ have the same sign (i.e., are either both nonnegative, or  both nonpositive).    Then $L_K(E) \cong L_K(F)$ as $K$-algebras.

\medskip

We note that, as of Fall 2013, it is not known whether the hypothesis regarding the germane  determinants can be eliminated from the statement of The Algebraic KP Theorem.   On the other hand, is the case that the determinant hypothesis {\it can} be eliminated in the analogous graph C$^*$-algebra result established by Kirchberg and Phillips.     Thus \cite[Corollary 2.7]{ALPS} is sometimes called the `Restricted' Algebraic Kirchberg Phillips  Theorem.

\medskip

Our goal for the remainder of this short note is to analyze the data required to invoke The Algebraic KP  Theorem  in the context of the collection of algebras $\{ L_K(C_n) \ | \ n\in \mathbb{N}\}$.  We start by noting that displays $(***)$ and $(**)$ immediately give

%{\bf Proof.}   The graph $C_n$ has the property that every vertex connects to every other vertex via some directed path, and that every cycle in $C_n$ has an exit in $C_n$.  (This last statement requires that $n\geq 3$.)    The second result follows from Enrique Pardo's direct argument, or by using properties of purely infinite simple algebras together with the description of purely infinite simple Leavitt path algebras and their $K_0$ groups.

\begin{proposition}\label{LKCnpis}
For each $n\geq 1$ the $K$-algebra $L_K(C_n)$ is unital purely infinite simple.  In particular, $M_{C_n}^* = (M_{C_n} \setminus \{[z]\},+)$ is a group. 
\end{proposition}

Referring to the explicit description of $M_{C_3}$ given above, it is easy to see that $$M_{C_3}^* = \{ [v_1], \  [v_2], \ [v_3], \ [v_1]+[v_2]+[v_3] \} \cong \Z/2\Z \times \Z/2\Z.$$
Using the previous computations, we see that $x = [v_1]+[v_2]+[v_3 ]$ is the identity element of the group $M_{C_3}^*$.  Indeed, we will see below that for any $n\in \mathbb{N}$,  the element $\sum_{i=1}^n [v_i]$ is the identity element in $M_{C_n}^*$.  (However, for an arbitrary graph $E$, the analogous element $\sum_{v\in E^0}[v]$ need not be the identity of $M_E^*$.)  

\medskip

Let $E$ be a finite directed graph for which $|E^0| = n$, and let $A_E$ denote the usual incidence matrix of $E$.    Let $B_E$ denote the matrix $I_n - A_E^t$. 
%(where $I_n$ is the identity $n \times n$ matrix, and $( \ )^t$ denotes the transpose of a matrix).    
We view $B_E$ both as a matrix, and as a linear transformation $B_E: \Z^n \rightarrow \Z^n$,  via left multiplication (viewing elements of $\Z^n$ as column vectors).   
In the situation where $L_K(E)$ is purely infinite simple, so that in particular $M_E^*$ is a group (necessarily isomorphic to $K_0(L_K(E))$), we have that  $$K_0(L_K(E)) \cong M_E^* \cong \Z^n / {\rm Im }(B_E) = {\rm Coker}(B_E).$$
(See \cite[Section 3]{AALP} for a complete discussion.)  Under this isomorphism we have  $[v_i ] \mapsto \vec{b_i} + {\rm Im }(B_E)$, where  $\vec{b_i}$ is the element of $\Z^n$ which is $1$ in the $i^{th}$ coordinate and $0$ elsewhere.

For any $n\times n$ matrix $T \in {\rm M}_n(\Z)$, we may view $T$ as a linear transformation from $\Z^n$ to $\Z^n$ in the usual way.  Then the finitely generated abelian group $\Z^n / {\rm Im }(T)$ may be described by analyzing the {\it Smith normal form} of $T$ (see e.g. \cite{SNFref}).     Specifically, if the Smith normal form of $T$ is the diagonal matrix $\rm{diag}(\alpha_1, \alpha_2, ..., \alpha_{n})$, then $\Z^n / {\rm Im }(T) \cong \Z / \alpha_1 \Z \oplus \Z / \alpha_2 \Z \oplus \cdots \oplus \Z / \alpha_n \Z$, where $\Z / 1\Z$ is interpreted as the trivial group $\{0\}$.   (Most computer software packages, e.g. {\it Mathematica} and {\it Scientific Notebook},  contain a built-in Smith normal form function.)    Using this method, where we let $T$ be the matrix $B_E = I_n - A_E^t$, we present here a description of the groups $M_{C_1}^*$ through $M_{C_{12}}^*$.  (Of course the description of some of these groups can be achieved using a more straightforward approach than the utilization of Smith normal form, as was done above for the group $M_{C_3}^*$.) 
\small 
$$M_{C_1}^* \cong \{0\}, \  M_{C_2}^* \cong \Z/3\Z,  \ M_{C_3}^* \cong \Z/2\Z \times \Z/2\Z,  \  M_{C_4}^* \cong \Z/3\Z, \  M_{C_5}^* \cong \{0\}, \ M_{C_6}^* \cong \Z \times \Z, $$
$$M_{C_7}^* \cong \{0\}, \   M_{C_8}^* \cong \Z/3\Z, \  M_{C_9}^* \cong \Z/2\Z \times \Z/2\Z, \  M_{C_{10}}^* \cong \Z/3\Z,    M_{C_{11}}^* \cong \{0\},   M_{C_{12}}^* \cong \Z \times \Z.$$
\normalsize

The displayed isomorphisms suggest a pattern, first noticed by A. Egri-Nagy and shared with the first author  in a private communication.

 %We denote the identity element of $M_{C_n}^*$ by $e$.    
 
 \medskip
 
 Although for an arbitrary graph $E$ there is no appropriate notion of ``the inverse of an element" in the semigroup $M_E^*$ (since this need not be a group), we use the standard minus sign notation to denote inverses in $M_E^*$ whenever this semigroup is actually a group.  In particular, it is appropriate to consider an element of the form $-[x]$ in  $M_{C_n}^*$ for any $[x] \in M_{C_n}^*$.

\begin{lemma}\label{mod6lemma}
For each $[v_i] $ in the group $M_{C_n}^*$ we have $[v_i ]= -[v_{i+3}]$.  Consequently, $[v_i ]= [v_{i+6}]$ for all $1 \leq i \leq n$ in $M_{C_n}^*$.
\end{lemma}

{\bf Proof.}   In $M_{C_n}^*$ we have $[v_{i+1}]= [v_{i}] + [v_{i+2}]$ and $[v_{i+2}]= [v_{i+1}] + [v_{i+3}]$.  Substituting yields  $[v_{i+1}] = [v_{i}] + [v_{i+1}] + [v_{i+3}].$   Since $M_{C_n}^*$ is a group, we can cancel $[v_{i+1}]$, and the result follows immediately.    \hfill $\Box$

%\end{enumerate}

%\begin{lemma}  

%\begin{enumerate}
%\item  In $M_6^*$, $\sum_{i=1}^6 v_i = 0$.
%\item  In each $M_n^*$ we have $\sum_{i=1}^{6k} = 0$ for any integer $k$ having $6k \leq n$.    Moreover, 
%\item  In $M_7^*$, $v_7 = 0$.
%\item  In $M_8^*$, $v_7 + v_8 = 0$.
%\item  In $M_9^*$, $v_7 + v_8 + v_9 = 0$.
%\item  In $M_{10}^*$, $v_7 + v_8 + v_9 + v_{10}= 0$.  (Indeed, $v_7 + v_8 = 0 = v_9 + v_{10}$.)
%\item  In $M_{11}^*$, $v_7 + v_8 + v_9 + v_{10} + v_{11} = 0$.  (Indeed, $v_7 = v_8  = v_9 = v_{10} = v_{11} = 0$.)
%\end{enumerate}
%\end{lemma}

\begin{proposition}\label{isomorphismsbetweengroups}
For $n, m\geq 1$, if $n \equiv m \ {\rm mod} \ 6$, then $M_{C_n}^* \cong M_{C_m}^*$.   
\end{proposition}

{\bf Proof.}    By the previously displayed isomorphisms, it suffices to show that if $n\geq 6$ and $n \equiv N {\rm mod} \ 6$ with $6\leq N \leq 11$, then $M_{C_n}^* \cong M_{C_N}^*$.   (We choose $N\geq 6$ in order to avoid some notational issues involving the interpretation of integers ${\rm mod}6$.) 
%establish the resultwe need establish the result in general for $n, n'\geq 6$, and then show that the cases $1\leq n \leq 5$ are valid as well.       For $n \geq 6$  let $N$ be the integer for which $n \equiv N {\rm mod} \ 6$, and $6 \leq N \leq 11$.  

%Recall that $B_j$ denotes the matrix  $ I_j - A_j^t$ for each positive integer $j$.  Using the previously displayed isomorphisms, it suffices to show that   $\Z^n / {\rm Im }(B_n) \cong \Z^N / {\rm Im }(B_N)$.   

%For each positive integer $j$ Let $F_t$ denote the free abelian monoid on generators $v_1, v_2, ..., v_t$, and   let $R_t$ denote the submonoid of $F_t$ generated by the expressions $v_i - v_{i-1} - v_{i+1}$ ($1\leq i \leq t$), where indices are interpreted mod $t$.   In particular, $M_t \cong F_t / R_t$, and $M_t^* \cong (F_t / R_t)^* $.   

For any graph $E$ having $|E^0| = t$ we denote by $\pi_t$  the canonical homomorphism of monoids $F_t \rightarrow F_t / R_t(E)$.    We note that since none of the generating relations which produce $R_t(E)$ are of the form $\sum_{i=1}^n m_iv_i = 0$ for $m_i \in \mathbb{N}$,  $\pi_t$ restricts to a homomorphism 
%(of semigroups, although the codomain is a group)
  $\pi_t^*: F_t^* \rightarrow (F_t / R_t(E))^*.$

  Define the semigroup homomorphism $\varphi: F_n^* \rightarrow F_N^*$ by setting $\varphi(v_i) = v_{i  {\rm mod} 6} $ for each generator $v_i$ ($1\leq i \leq n$) of $F_n^*$, and extending linearly.   Let $\psi: F_n^* \rightarrow M_{C_N}^*$ be the composition $\pi_N^* \circ \varphi$, so that, in particular, $\psi(v_i) = [v_{i {\rm mod}6}]$ for $1\leq i \leq n$.  
To show that $\psi$ factors to a homomorphism from $M_{C_n}^*$ to $M_{C_N}^*$, we need only show that $\psi$ takes each of the relations $v_i =  v_{i-1} + v_{i+1}$ ($1\leq i \leq n$, interpreted ${\rm mod}n$) in $R_n(E)$ to a relation valid in $M_{C_N}^*$;  in other words, it suffices to show that  $\psi(v_i )= \psi(  v_{i-1} )+ \psi(v_{i+1})$ in $M_{C_N}^*$  We consider five cases.  The point here is that we must understand the given relations with two things in mind:   the interpretation of the integers $k {\rm mod} 6$ which arise in the definition of $\varphi$ as integers between $1$ and $6$ (inclusive), as well as the interpretation of the subscripts in the expressions in $M_{C_N}^*$ as integers ${\rm mod}N$.  

\smallskip

Case 1a:   $1 \leq i-1$ and $i+1 \leq n$ and $i-1 \equiv 1, 2, 3, $ or $4, \ {\rm mod} \ 6$.   Then 
$$\psi(v_i ) =  [v_{i  {\rm mod} 6}] =   [v_{i-1  {\rm mod} 6}]  + [v_{i +1  {\rm mod} 6}] = \psi(v_{i-1} ) + \psi (v_{i+1}).$$

\smallskip

Case 1b:  $1 \leq i-1$ and $i+1 \leq n$ and   $i-1 \equiv 5  {\rm mod} 6$; so   $i \equiv 6  {\rm mod} 6$ and $i +1  \equiv 1  {\rm mod} 6$. Then using Lemma \ref{mod6lemma} we have
$$\psi(v_i ) =  [v_{i  {\rm mod} 6}] =  [v_6]  = [v_5] + [v_7] = [v_5] + [v_1]  =   [v_{i-1  {\rm mod} 6}]  + [v_{i +1  {\rm mod} 6}]$$
$$  \hspace{-2.45in} = \psi(v_{i-1} ) + \psi (v_{i+1}).$$

\medskip

Case 1c:  $1 \leq i-1$ and $i+1 \leq n$ and   $i-1 \equiv 6  {\rm mod} 6$; so   $i \equiv 1  {\rm mod} 6$ and $i +1  \equiv 2  {\rm mod} 6$. Then using Lemma \ref{mod6lemma} we have
$$\psi(v_i ) =  [v_{i  {\rm mod} 6}] =  [v_1]  = [v_7] =  [v_6] + [v_8] =  [v_6] + [v_2] $$
$$ 
%\hspace{-2.5in}
\ \ \ \ \ \    =  [v_{i-1  {\rm mod} 6}]  + [v_{i +1  {\rm mod} 6}] 
= \psi(v_{i-1} ) + \psi (v_{i+1}).$$

\medskip

Case 2:   $i=1$.  So $v_{i-1} = v_n$ in $F_n$ by definition.   Then using that $N {\rm mod}6 \equiv n {\rm mod}6$ and Lemma \ref{mod6lemma}, we get
$$\psi(v_1) =  [v_1] =  [v_N] + [v_2] = [v_{N {\rm mod}6}] + [v_2]  =  [v_{n {\rm mod}6}] + [v_2]$$
$$  \hspace{-.1in}  =[v_{i-1  {\rm mod} 6}]  + [v_{i +1  {\rm mod} 6}] = \psi(v_{i-1} ) + \psi (v_{i+1}).$$

\medskip

Case 3:   $i=n$.  So  $v_{i+1} = v_1$ in $F_n$ by definition.  Then using that $N {\rm mod}6 \equiv n {\rm mod}6$ and Lemma \ref{mod6lemma}, we get
$$\psi(v_n) =  [v_{n {\rm mod} 6}] = [v_{N {\rm mod} 6}]  = [v_N] = [v_{N-1}] + [v_1] = [v_{N-1 {\rm mod}6}] + [v_1] 
$$
$$ \ \ \ \ \ \ \ \ \ \  =  [v_{n -1 {\rm mod}6}] + [v_1] 
 =[v_{i-1  {\rm mod} 6}]  + [v_{i +1  {\rm mod} 6}] = \psi(v_{i-1} ) + \psi (v_{i+1}).$$

\medskip

%$$\pi_N(v_{i  {\rm mod} 6} - v_{i-1  {\rm mod} 6} - v_{i +1  {\rm mod} 6}) =  \pi_N( v_{6} - v_{5 } - v_{7}$$
 
 % We claim that ${\rm Im}(B_n) \subseteq {\rm Ker}(\psi)$.   So suppose that $\vec{x} \in \Z^n$ has $\vec{x} = \vec{y} - A_n^t\vec{y}$ for some $\vec{y} \in \Z^n$.   

Thus $\psi$ preserves the relations which generate $R_n(E)$, and so  $\psi$ extends to a group homomorphism $\overline{\psi}: M_{C_n}^* \rightarrow M_{C_N}^*$.

In a completely analogous manner, for $n \equiv N {\rm mod}6$ with $6\leq N \leq 11$ there exists a group homomorphism $\overline{\tau}: M_{C_N}^* \rightarrow M_{C_n}^*$ for which $\overline{\tau}([v_i]) = [v_{i {\rm mod}6}]$.  It   is then clear that $\overline{\tau}$ and $\overline {\psi}$ are inverses, thus establishing the result for $n,m \geq 6$.

The observation made prior to the Proposition shows that the cases $n=1, 2,3,4,5$ satisfy the statement as well.  \hfill $\Box$

\medskip

As a consequence of Proposition \ref{isomorphismsbetweengroups}, and using the previously displayed computations, we see that there are, up to isomorphism, only four groups represented by the collection $\{M_{C_n}^*  \ | \ n\in \mathbb{N}\}$, as follows.   
%$$M_{C_3}^* \cong \Z/2\Z \times \Z/2\Z,    \ \ \ M_{C_4}^* \cong \Z/3\Z,  \ \ \ M_{C_5}^* \cong \{0\},  \ \ \ M_{C_6}^* \cong \Z \times \Z. $$
%$$M_{C_1}^* \cong \{0\},  \ \ \ M_{C_2}^* \cong \Z/3\Z,  \ \ \ M_{C_3}^* \cong \Z/2\Z \times \Z/2\Z,    \ \ \ M_{C_4}^* \cong \Z/3\Z,  \ \ \ M_{C_5}^* \cong \{0\},  \ \ \ M_{C_6}^* \cong \Z \times \Z, $$

\begin{corollary}\label{fourisoclassescorollary}  The following is a complete description of the isomorphism classes of the groups $M_{C_n}^*$ for $n\in \mathbb{N}$. 
\begin{enumerate}
\item $M_{C_n}^* \cong \{0\}$ in case $n \equiv 1 {\rm mod} 6$ or $n \equiv 5 {\rm mod} 6$.
 \item $M_{C_n}^* \cong \Z / 3\Z$ in case $n \equiv 2 {\rm mod} 6$ or $n \equiv 4 {\rm mod} 6$.
\item $M_{C_n}^* \cong \Z / 2\Z \times \Z / 2\Z$ in case $n \equiv 3 {\rm mod} 6$.
\item $M_{C_n}^* \cong \Z  \times \Z $ in case $n \equiv 6 {\rm mod} 6$.
\end{enumerate}
\end{corollary}

\smallskip

With the first piece of the analysis now in place, we turn our attention to describing the element $[L_K(C_n)]$ of the group $K_0(L_K(C_n))$; by $(*)$, this amounts to describing the element $\sum_{i=1}^n [v_i]$ in the group $M_{C_n}^* $.   
\begin{lemma}\label{sumverticesisidentity}
For each $n\in \mathbb{N}$, $\sum_{i=1}^n [v_i ]$ is the identity element of the group  $M_{C_n}^*$.
\end{lemma}

{\bf Proof}:   Let $x$ denote the element $\sum_{i=1}^n [v_i ]$ of $M_{C_n}^*$.  Using the defining relations $R_n(C_n)$, we have  
$$ x = \sum_{i=1}^n [v_i ] = \sum_{i=1}^n ([v_{i-1} ] + [v_{i+1} ]) = \sum_{i=1}^n [v_{i-1} ] + \sum_{i=1}^n [v_{i+1} ] = x + x,$$
as we interpret the indices of the generating elements modulo $n$ in $M_{C_n}^*$.  But the equation $x+x=x$ in a group yields immediately that $x$ is the identity element.   \hfill $\Box$

\medskip

Suppose that $m,n$ are integers for which $M_{C_n}^* \cong M_{C_m}^*$ (see  Corollary \ref{fourisoclassescorollary}).   Then trivially such an isomorphism must send the element $\sum_{i=1}^n [v_i]$ of $M_{C_n}^*$ to the element  $\sum_{i=1}^m [v_i]$ of $M_{C_m}^*$, as  Lemma  \ref{sumverticesisidentity}  shows that each of these expressions is the identity element in the respective group.   

The final piece of the hypotheses in the Algebraic KP Theorem involves determinants of appropriate matrices, which we analyze in the next result.

\begin{proposition}\label{nonpositivedet}
For each $n\in \mathbb{N}$, ${\rm det} (I_n - A_{C_n}^t) \leq 0$.   
\end{proposition} 
%%%%%%% BEN STUFF
{\bf Proof}:   An $n\times n$ matrix $B = (b_{i,j})$ is {\it circulant} in case we have $b_{i+1, j+1} = b_{i,j}$ for  $1\leq i \leq n-1$, $1\leq j \leq n-1$;  $b_{1,j+1} = b_{n,j}$ for $j+1 \leq n$; and $b_{i+1,1} = b_{i,n}$ for $i+1\leq n$.   Less formally, $B$ is circulant in case each subsequent row of $B$ is obtained from the previous row by moving each entry of the previous row one place to the right, and moving the last entry of the previous row to the first position of the subsequent row.   (The last row gets moved to the first row in this way as well.)   If $B$ is circulant, then there is a formula (derived from an analysis of the eigenvectors of $B$) which expresses  ${\rm det}(B) $ as the following product:  
$$ \det(B) = \prod\limits_{j=0}^{n-1}(b_1+b_{2}\omega_j+b_{3}\omega_j^2+\dotsb + b_n\omega_j^{n-1})$$ where  $(b_1 \ b_2 \ b_3 \ \cdots \ b_n)$ is the first row of $B$, and $\omega_j = e^{\frac{2\pi i j }{n}}$ is an  $n$th root of unity in $\mathbb{C}$.   (See e.g. \cite{circulantref} for a  description of some of the many places in which circulant matrices arise.) 
%, and $\omega_j=\exp (\frac{2\pi ij}{n})$.
%and $i=\sqrt{-1}$, the imaginary unit.\\[4mm]

In the case of the  Cayley graph $C_n$ (for $n\geq 3$),  the matrix $B =  I_n-A_{C_n}^t$ has $b_1=1$, $b_2=b_{n}=-1$, and $b_i=0$ for $i= 3,  4, \dotsc n-1$.  Using that $e^{i\theta} =\cos\theta+i\sin\theta$ together with the displayed equation,  we get:
\begin{equation*}
\begin{split}
\det(I_n-A_{C_n}^t) 
%& = \prod\limits_{j=0}^{n-1}(b_1+c_{n-1}\omega_j+c_{n-2}\omega_j^2+\dotsb + c_1\omega_j^{n-1})\\
& = \ \prod\limits_{j=0}^{n-1}(1-\omega_j-\omega_j^{n-1}) \ = \ 
 \prod\limits_{j=0}^{n-1}(1-e^{\frac{2\pi ij}{n}}-e^{\frac{2\pi ij(n-1)}{n}})\\
& = \  \prod\limits_{j=0}^{n-1}(1-\cos\frac{2\pi j}{n}- i\sin\frac{2\pi j}{n}-\cos\frac{2\pi j(n-1)}{n}-i\sin\frac{2\pi j(n-1)}{n})\\
& = \ \prod\limits_{j=0}^{n-1}(1-2\cos\frac{2\pi}{n}j),\\
\end{split}
\end{equation*}
\noindent
with the final equality coming as a direct result of the basic trigonometry facts that, for any integer $j$,    $\cos\frac{2\pi j (n-1)}{n}=\cos\frac{2\pi j}{n}$  and  $\sin\frac{2\pi j(n-1)}{n}=-\sin\frac{2\pi j}{n}$.  

  When $j=0$ we have  $1-2\cos(\frac{2\pi }{n} j) = 1-2\cdot \cos 0 =-1<0$.   In case $n$ is even,  when $j=\frac{n}{2}$ we have $1-2\cos(\frac{2\pi }{n} j) =1-2\cos\pi=3>0$. Furthermore, since 
  %$\cos\frac{2\pi}{n}k=\cos\frac{2\pi}{n}(n-k),$
   $1-2\cos\frac{2\pi}{n}j=1-2\cos\frac{2\pi}{n}(n-j)$, we see that  $(1-2\cos\frac{2\pi}{n}j ) ( 1-2\cos\frac{2\pi}{n}(n-j) )
   %= (1-2\cos\frac{2\pi}{n}j )^2 
  \geq 0$.   This yields the result.   
  %Thus when $n$ is odd, the terms in the product corresponding to $j \neq 0$  occur in pairs, specifically, that   Since the $j=0$ term is $-1$, we have    $\Pi_{j=0}^{n-1}(1-2\cos\frac{2\pi}{n}j)\leq 0$ for $n$ odd.  On the other hand, when $n$ is even, then  The other terms of the product, again with the exception of $j=0$, occur in pairs as in the previous case.  So again, we have $\Pi_{j=0}^{n-1}(1-2\cos\frac{2\pi}{n}j)\leq 0$. 
     \hfill $\Box$

\medskip

We note as a consequence of the previous analysis that ${\rm det}(I_n - A_{C_n}^t) = 0$ precisely when one of the factors $1-2\cos\frac{2\pi}{n}j$ ($0\leq j \leq n-1$)  equals $0$.   This can easily be shown to happen precisely when $n$ is a multiple of $6$.   This information is consistent with the observation that the only values of $n$ for which the group $M_{C_n}^*$ is infinite are multiples of $6$.

%%%%%%%%%

Of the four groups which arise up to isomorphism as a group of the form $M_{C_n}^*$ (for $n\geq 1$), we see that two are cyclic:   $\{0\}$ and $\Z/3\Z$.    Purely infinite simple unital Leavitt path algebras $L_K(E)$ whose corresponding $K_0$ groups are cyclic and for which ${\rm det}(I_{|E^0|}  - A_E^t) \leq 0$ are relatively well-understood, and arise from the classical {\it Leavitt algebras} $L_K(1,n)$, as follows.    For any integer $n\geq 2$, $L_K(1,n)$  is the free associative $K$-algebra  in $2n$ generators $x_1, x_2, ..., x_n, y_1, y_2, ..., y_n$, subject to the relations
$$ y_i x_j = \delta_{i,j}1_K \ \ \mbox{and} \ \ \sum_{i=1}^n x_i y_i = 1_K.$$
These algebras were first defined and investigated in \cite{L}, and formed the motivating examples for the more general notion of  Leavitt path algebra.  It is easy to see  that for $n\geq 2$,  if $R_n$ is the graph having one vertex and $n$ loops (the ``rose with $n$ petals" graph), then $L_K(R_n) \cong L_K(1,n)$.    It is clear from $(***)$ that each $L_K(R_n)$ is purely infinite simple; it is straightforward from $(**)$ that  $K_0(L_K(R_n)) \cong M_{R_n}^*$ is the cyclic group $\Z / (n-1)\Z$ of order $n-1$, where the regular module $[L_K(R_n)]$ in $K_0(L_K(R_n))$ corresponds to $ 1$ in $\Z / (n-1)\Z$.   

Now let $d \geq 2$, and consider the graph $R_n^d$ having two vertices $v_1, v_2$; $d-1$ edges from $v_1$ to $v_2$; and $n$ loops at $v_2$:  
$$  R_n^d \ = \   \xymatrix{
  \bullet^{v_1}  \ar[r]^{(d-1)}
 & \bullet^{v_2}  \ar@(ur,dr)[]^{(n)}  } \ \ \ \  \
$$
   It is shown in \cite{AALP} that  the matrix algebra ${\rm M}_d(L_K(1,n))$ is isomorphic to $L_K(R_n^d)$.  By standard Morita equivalence theory, 
   %this implies that 
   $K_0({\rm M}_d(L_K(1,n))) \cong K_0(L_K(1,n)).$
   % \cong \Z / (n-1)\Z$. 
   Moreover, the element $[{\rm M}_d(L_K(1,n))]$ of $K_0({\rm M}_d(L_K(1,n)))$  corresponds to the element $d$ in $ \Z / (n-1)\Z$.   In particular,  the element $[{\rm M}_{n-1}(L_K(1,n))]$ of $K_0({\rm M}_{n-1}(L_K(1,n)))$  corresponds to $n-1 \equiv 0$ in $  \Z / (n-1)\Z$.    Finally, an easy computation yields that ${\rm det} (I_2 - A_{R_n^d}^t) = -(n-1) <0 $ for all $n,d$.    
Therefore, by invoking the Algebraic KP Theorem, the previous discussion immediately yields the following.  

\begin{proposition}\label{isotomatrixoverLeavitt}
Suppose $E$ is a graph for which $L_K(E)$ is unital purely infinite simple.  Suppose that $M_E^*$ is isomorphic to the cyclic group $ \Z / (n-1)\Z$, via an isomorphism which takes the element $ \sum_{v\in E^0}[v]$ of $M_E^*$  to the element $d$ of  $ \Z / (n-1)\Z$.  Finally, suppose that ${\rm det} (I_{|E^0|} - A_E^t)$ is negative.    Then
$L_K(E) \cong {\rm M}_d(L_K(1,n)).$
\end{proposition}

We now have all the ingredients in place to achieve our main result.  

\begin{theorem}\label{MainTheorem}
For each $n \geq 1$ let $C_n$ denote the Cayley graph  corresponding to the cyclic group $\Z/n\Z$ (with respect to the subset $\{1,n-1\}$) as described previously.   Then up to isomorphism the collection of Leavitt path algebras $\{L_K(C_n) \ | \ n\in \mathbb{N}\}$ is completely described by the following four pairwise non-isomorphic classes of $K$-algebras. 

 \begin{enumerate}

\item $L_K(C_n) \cong L_K(C_m)$  in case $ m \equiv 1 $ or $5  \ {\rm mod} 6$  and $n \equiv 1 $ or $5 \  {\rm mod} 6$.  In this case, these algebras are isomorphic to $L_K(1,2)$.

\item $L_K(C_n) \cong L_K(C_m)$  in case $ m \equiv 2 $ or $4  \ {\rm mod} 6$  and $n \equiv 2 $ or $4 \  {\rm mod} 6$.  In this case, these algebras are isomorphic to ${\rm M}_3(L_K(1,4))$.

\item $L_K(C_n) \cong L_K(C_m)$  in case $ m, n  \equiv 3 \ {\rm mod} 6$.  

\item $L_K(C_n) \cong L_K(C_m)$  in case $ m, n  \equiv 6 \ {\rm mod} 6$.

\end{enumerate}

\end{theorem}

{\bf Proof}:   We seek to invoke the Algebraic KP Theorem.   By Proposition \ref{LKCnpis}, $L_K(C_n)$ is purely infinite simple for each $n\in \mathbb{N}$.   For any of the four indicated cases,  we choose a pair of integers $m,n$.   By Corollary \ref{fourisoclassescorollary} and $(**)$ we have $K_0(L_K(C_n)) \cong K_0(L_K(C_m))$, and, by the observation made subsequent to Lemma \ref{sumverticesisidentity} together with $(*)$,  this isomorphism necessarily takes $[L_K(C_n)]$ to  $[L_K(C_m)]$.   By Proposition \ref{nonpositivedet},  ${\rm det}(I_n - A_{C_n}^t)$ and ${\rm det}(I_m - A_{C_m}^t)$ are both nonpositive.    The Algebraic KP Theorem now gives the bulk of the  result.  The two extra statements in parts (i) and (ii) follow directly from statements (i) and (ii) of Corollary \ref{fourisoclassescorollary} together with Proposition \ref {isotomatrixoverLeavitt}.   \hfill $\Box$

\medskip

In \cite{AAP3} a description is given of the Leavitt path algebras associated to additional  collections  of Cayley-type graphs.

%\section*{This is an unnumbered first-level section head}
%This is an example of an unnumbered first-level heading.

%\specialsection*{This is a Special Section Head}
%This is an example of a special section head%
%%%%%%%%%%%%%%%%%%%%%%%%%%%%%%%%%%%%%%%%%%%%%%%%%%%%%%%%%%%%%%%%%%%%%%%%
%\footnote{Here is an example of a footnote. Notice that this footnote
%text is running on so that it can stand as an example of how a footnote
%with separate paragraphs should be written.
%\par
%And here is the beginning of the second paragraph.}%
%%%%%%%%%%%%%%%%%%%%%%%%%%%%%%%%%%%%%%%%%%%%%%%%%%%%%%%%%%%%%%%%%%%%%%%%
.

%\section{This is a numbered first-level section head}
%This is an example of a numbered first-level heading.

%\subsection{This is a numbered second-level section head}
%This is an example of a numbered second-level heading.

%\subsection*{This is an unnumbered second-level section head}
%This is an example of an unnumbered second-level heading.

%\subsubsection{This is a numbered third-level section head}
%This is an example of a numbered third-level heading.

%\subsubsection*{This is an unnumbered third-level section head}
%This is an example of an unnumbered third-level heading.

\bibliographystyle{amsalpha}

\end{document}